\documentclass[amsthm]{elsart}

\usepackage{yjsco}
\usepackage{natbib}
\usepackage{amssymb}
\usepackage{amsmath}
\usepackage{algorithm}
\usepackage[noend]{algpseudocode}

\begin{document}
\begin{frontmatter}

\title{The Termination of Algorithms for Computing Gr\"{o}bner Bases}

\thanks{This work was supported by the National Natural Science Foundation of
China (No. 61173151, No. 61173152) and the 111 Project (No. B08038).}

\author[ISN]{Senshan Pan\corauthref{cor}},
%\address{State Key Laboratory of Integrated Service Networks, Xidian University, Xi'an 710071, China}
\corauth[cor]{Corresponding author.}
\ead{pansenshan@gmail.com}
%\ead[url]{URL 1}
\author[ISN]{Yupu Hu} and
%\address{State Key Laboratory of Integrated Service Networks, Xidian University, Xi'an 710071, China}
\ead{yphu@mail.xidian.edu.cn}
%\ead[url]{URL 2}
\author[ISN]{Baocang Wang}
\ead{bcwang79@yahoo.com.cn}

\address[ISN]{State Key Laboratory of Integrated Service Networks, Xidian University, Xi'an 710071, China}

\begin{abstract}
The F5 algorithm \citet{Fau02} is generally believed as one of
the fastest algorithms for computing Gr\"{o}bner bases. However, its termination
problem is still unclear. Recently, an algorithm GVW \citet{GVW}
and its variant GVWHS \citet{Volny} have been proposed, and their efficiency
are comparable to the F5 algorithm. In the paper, we clarify the concept of an admissible
module order. For the first time, the connection between the reducible
and rewritable check is discussed here. We show that the top-reduced S-Gr\"{o}bner basis must be finite if the admissible monomial order and the admissible module order are compatible. Compared with \citet{Volny}, this paper presents a complete proof of the termination and correctness of the GVWHS algorithm. What
is more, it can be seen that the GVWHS is in fact an F5-like algorithm.
Different from the GVWHS algorithm, the F5B algorithm may generate redundant sig-polynomials. Taking into account this situation, we prove the termination and correctness of
the F5B algorithm. And we notice that the original F5 algorithm in
\citet{Fau02} slightly differs from the F5B algorithm in the insertion strategy on which the F5-rewritten criterion is based.
Exploring the potential ordering of sig-polynomials computed by the original F5 algorithm, we propose an F5GEN algorithm with a generalized insertion strategy, and prove the termination and correctness of it. Therefore, we have a positive answer to the long standing problem of proving the termination of the original F5 algorithm.
\end{abstract}

\begin{keyword}
Termination \sep Gr\"{o}bner basis \sep GVWHS \sep F5
\end{keyword}

\end{frontmatter}

%\section{Main Result}
%
%We recall the works reported in \citetp{Key1}.
%......
%We recall that \citett{Key1} proved the result.
%
%\begin{defn}[Nice Notion]
%We say a polynomial is \emph{nice} if ....  \qed
%\end{defn}
%
%We state the main result.
%
%\begin{thm}[Someone 2005]
%All nice polynomials are also pretty....  \qed
%\end{thm}
%
%\begin{pf}
%The proof is easy...
%\end{pf}
%
%\begin{ack}
%The authors would like to thank ....
%\end{ack}
%
%\bibliographystyle{elsart-harv}
%\bibliography{mybibfile}

\section{Introduction}
\label{intro}

In cryptography, the cipher of a cryptosystem sometimes can be transformed into a system of equations.
Solving a set of multivariate polynomial equations (nonlinear and randomly chosen) over a finite field is an NP-hard problem \citet{Garey79}. Based on which, %a stream cipher with provable security named QUAD \citet{Berbain} was proposed. Recently,
Albrecht \textit{et al.} \citet{Albrecht11} constructed a Polly-Cracker-style cryptosystem. However, in much more cases a designer
has to embed some kind of trapdoor function to enable efficient decryption and signing. Although the structure of the cipher is hidden, the equations are so special that one can exploit them via Gr\"{o}bner basis based techniques to attack the cryptosystem.

In 1965 Buchberger's \citet{Buch65} thesis he described the appropriate
framework for the study of polynomial ideals, with the introduction
of Gr\"{o}bner bases. Since then, Gr\"{o}bner basis has become a fundamental
tool of computational algebra and it has found countless applications
in coding theory, cryptography and even directions of Physics, Biology
and other sciences.

Although Buchberger presented several improvements to his algorithm
for computing Gr\"{o}bner bases in \citet{Buch79}, the efficiency
is not so good. Recent years have seen a surge in the number of algorithms
in computer algebra research, but efficient ones are few.
Faug\`{e}re \citet{Fau02} proposed the idea of signatures and utilized
two powerful criteria to avoid useless computation in the F5 algorithm.
Faug\`{e}re and Joux broke the first
Hidden Field Equation (HFE) Cryptosystem Challenge (80 bits) by using
the F5 algorithm in \citet{Fau03}. The proof of the termination in \citet{Fau02} was labeled as a conjecture in \citet{Stegers}. However, Gash \citet{Gash} pointed out that there
exists an error in the proof of the termination of the F5 algorithm, and he proposed another conjecture for it. It will be shown in this paper that the conjecture is still wrong. In
\citet{Arri11}, a simpler algorithm was constructed to
prove the termination, but the proof unfortunately has flaws due to
the abuse of the monomial order and the module order mentioned
in this paper. Though the F5 algorithm seems to terminate for any
polynomial ideals, the proof of it has been admitted as an open problem
in \citet{SunISSAC11}, \citet{EderISSAC11}, \citet{EderCCA11}.
Recently, signature-based algorithms like the GVW algorithm and
its variant the GVWHS algorithm are proposed in \citet{GVW}, \citet{Volny}.
The algorithms are claimed to terminate if the monomial order and
the module order are ``compatible'', but readers can hardly find a direct proof. The relation between the reducible and rewritable check, which was not considered before, is studied in the paper, and the finiteness
of the top-reduced S-Gr\"{o}bner basis for a polynomial ideal is proved if the ``compatible'' property is satisfied.
Then we give a complete proof of the termination of the GVWHS algorithm. Besides, through reformulation, the
GVWHS algorithm can be seen as an F5-like algorithm
(with a different insertion strategy). Though the F5B algorithm (F5 algorithm in Buchberger's style) may generate redundant sig-polynomials, by analyzing the similarity with the GVWHS algorithm, we prove the
termination of the F5B algorithm. Moreover, the termination of the F5GEN algorithm (F5 algorithm with a generalized insertion strategy) is also proved later on. Moreover, by employing an appropriate insertion strategy for the F5GEN algorithm , the proof of the correctness and termination of the original F5 algorithm is self-evident.
%At last, we make a comparison in time and

The paper is organized as follows. We start by settling basic
notations in Section \ref{pre}. In Section \ref{module order}, we present a new definition of the admissible
module order. Then two admissible orders and their connection are described in Section
\ref{sig-polynomial} and the top-reduced S-Gr\"{o}bner basis for a polynomial ideal is proved to be finite. Based on this finiteness, we propose a new proof of the termination of the GVWHS algorithm in Section \ref{the GVWHS} and point out that the GVWHS algorithm
is a variant of F5 algorithm by introducing the intermediate F5G algorithm (F5 algorithm in GVWHS's style). In Section \ref{the F5B}, a simpler version of the
F5B algorithm in \citet{SunAMSS09} is presented and proved. Considering the different insertion strategy between the
F5 algorithm and the F5B algorithm in this paper, we prove the correctness and termination
of the F5GEN algorithm in Section \ref{the F5}.

\section{Preliminaries}
\label{pre}

Let $r$ be a (binary) relation on a set $\mathbf{M}$, one may associate the \textbf{strict
part} $r_{s}=r\setminus r^{-1}$, and let $\mathbf{N}\subseteq \mathbf{M}$. Then an
element $a$ of $\mathbf{N}$ is called \textbf{$r$-minimal} in $\mathbf{N}$ if there
is no $b\in \mathbf{N}$ with $b\, r_{s}\, a$. A \textbf{strictly descending $r$-chain} in $\mathbf{M}$
is an infinite sequence $\{a_{n}\}_{n\in \mathbf{N}}$ of elements of $\mathbf{M}$
such that $a_{n+1}\, r_{s}\, a_{n}$ for all $n\in \mathbf{N}$. If there is another relation $t$ satisfying $r\subseteq t$, then $t$ is called an \textbf{extension} of r.
The relation
$r$ is called \textbf{well-founded} if every non-empty subset $\mathbf{N}$
of $\mathbf{M}$ has an $r$-minimal element, $r$ is a \textbf{well-order}
on $\mathbf{M}$ if $r$ is a well-founded linear order on $\mathbf{M}$. For more concepts
not presented here, refer to \citet{Becker}.

Let $R=k[x_{1},\ldots,x_{n}]$ be the polynomial ring over the field
$k$ with $n$ variables. We define an admissible order on the monoid
$\mathcal{M}=\{\Pi_{i=1}^{n}x_{i}^{a_{i}}\,|\, a_{i}\in \mathbb{N}\}$.

\begin{defn} An \textbf{admissible monomial order} $\le_{m}$ is a linear order on $\mathcal{M}$ that
satisfies the following conditions.
\begin{enumerate}
\item  $1\le_{m}m$ for all $m\in\mathcal{M}$.
\item  $m_{1}\le_{m}m_{2}$ implies $m_{1}\cdot s\le_{m}m_{2}\cdot s$
for all $s,m_{1},m_{2}\in\mathcal{M}$.
\end{enumerate}
\end{defn}

It can be seen that the admissible order $\leq_{m}$ is a well-order
on $\mathcal{M}$. Sometimes we write $=$ for $=_m$ for brevity. For any $p\in R$, without confusion, we denote
the leading monomial of $p$ by $lm(p)$, the leading coefficient by $lc(p)$,
and the leading term by $lt(p)=lc(p)lm(p)$ with respect to the order $\le_{m}$.

Let $\mathcal{I}$ be the ideal generated by the set $\mathbf{F}=\{f_{1},\ldots,f_{d}\}\in R$, that is,
\begin{equation*}
\mathcal{I}=<f_{1},\ldots,f_{d}>=\{p_{1}f_{1}+\ldots+p_{d}f_{d}\,|\, p_{1},\ldots,p_{d}\in R\}.
\end{equation*}
Consider the following $R$-submodule of $R^{d}\times R$:
\begin{equation*}
\mathcal{SP}=\{(\mathbf{u},p)\in R^{d}\times R\,|\,\mathbf{u}\cdot\mathbf{f}=p\},
\end{equation*}
where $\mathbf{f}=(f_{1},\ldots,f_{d})\in R^{d}$, and $\mathbf{e}_{i}$
is $i$-th unit vector of $R^{d}$ such that the free $R$-module
$R^{d}$ is generated by the set $\Sigma=\{\mathbf{e}_{1},\ldots,\mathbf{e}_{d}\}$.
The element $sp$ in $\mathcal{SP}$ we call a \textbf{sig-polynomial}.
A subset $\mathbf{Syz}=\{(\mathbf{u},0)\in\mathcal{SP}\}$ is defined
the \textbf{syzygy submodule} for $\mathcal{SP}$, and $\mathbf{NSP}=\mathcal{SP}\setminus \mathbf{Syz}$
is called the set of \textbf{non-syzygy sig-polynomials}. Let $(\mathbf{u}_1,p_1)$ and $(\mathbf{u}_2,p_2)$ be two non-syzygy sig-polynomials in $\mathcal{SP}$. A syzygy $(p_2\mathbf{u}_1-p_1\mathbf{u}_2,0)$ is called a \textbf{principal syzygy}.

\section{The admissible module order}
\label{module order}

Below is a fundamental tool for a clearer understanding
of termination of algorithms for computing Gr\"{o}bner bases,

\begin{defn} Let $\preceq$ be a quasi-order on $\mathbf{M}$ and let
$\mathbf{N}\subseteq\mathbf{M}$. Then a subset $\mathbf{B}$ of $\mathbf{N}$
is called a \textbf{Dickson basis}, or simply basis of $\mathbf{N}$
w.r.t. $\preceq$ if for every $a\in\mathbf{N}$ there exists some
$b\in\mathbf{B}$ with $b\preceq a$. We say that $\preceq$ has the
\textbf{Dickson property}, or is a \textbf{Dickson quasi-order}, if
every subset $\mathbf{N}$ of $\mathbf{M}$ has a finite basis w.r.t. $\preceq$. \end{defn}

If $\preceq$ is a (Dickson) quasi-order on $\mathbf{M}$, then we
call $(\mathbf{M},\preceq)$ a (Dickson) quasi-ordered set. Let now
$(\mathbf{M},\preceq)$ and $(\mathbf{N},\preceq)$ be quasi-ordered
sets, then a quasi-order $\preceq'$ on Cartesian product $\mathbf{M}\times\mathbf{N}$
is defined as follows:
\begin{equation*}
(a,b)\preceq'(c,d)\ \Leftrightarrow\ a\preceq b\ and\ c\preceq d,
\end{equation*}
for all $(a,b)$, $(c,d)$ $\in$ $\mathbf{M}\times\mathbf{N}$.
The \textbf{direct product} of the quasi-order sets $(\mathbf{M},\preceq)$
and $(\mathbf{N},\preceq)$ is denoted by $(\mathbf{M}\times\mathbf{N},\preceq')$.
The Dickson property can be derived as follows.

\begin{lem}\label{direct product}\citet{Becker}
Let $(\mathbf{M},\preceq)$ and $(\mathbf{N},\preceq)$ be Dickson
quasi-ordered sets, and let $(\mathbf{M}\times\mathbf{N},\preceq')$
be their direct product. Then $(\mathbf{M}\times\mathbf{N},\preceq')$
is a Dickson quasi-ordered set. \end{lem}

The immediate corollary is that $(\mathbb{N}^{n},\le')$, the direct
product of $n$ copies of the natural numbers $(\mathbb{N},\le)$
with their natural ordering is a Dickson partially ordered set. This
is Dickson's lemma, and another version of which is given below by an isomorphism.

\begin{lem}[Dickson's lemma] \label{Dickson's lemma}\citet{Becker}
The divisibility relation $|$ on $\mathcal{M}$ is a Dickson partial
order on $\mathcal{M}$. More explicitly, every non-empty subset $\mathbf{S}$
of $\mathcal{M}$ has a finite subset $\mathbf{B}$ such that for
all $s\in\mathbf{S}$, there exists $t\in\mathbf{B}$ with $t\,|\, s$.
\end{lem}

Let $\mathcal{M}_{d}=\{m\mathbf{e}_{i}\,|\, m\in\mathcal{M},i\in\{1,\ldots,d\}\}$
be the $\mathcal{M}$-monomodule of $R^{d}$. The definition of the
divisibility relation $|'$ on $\mathcal{M}_{d}$ is
\begin{equation*}
m_{1}\mathbf{e}_{i}\,|'\, m_{2}\mathbf{e}_{j}\leftrightarrow m_{1}|m_{2}\ and\ i=j\in\{1,\ldots,d\}.
\end{equation*}
 By an abuse of notation, we still denote $|$ instead of $|'$.
Since $(\mathcal{M},|)$ is a Dickson partial ordered set, by decomposing $\mathcal{M}_d$ into $\cup\mathcal{M}\mathbf{e}_i$, $(\mathcal{M}_{d}, |)$ is also a Dickson partial ordered set.
On $\mathcal{M}_{d}$, we will define the admissible order similarly.

\begin{defn} An \textbf{admissible module order} $\le_{s}$ is a linear
order on $\mathcal{M}_{d}$ that satisfies the following conditions.
\begin{enumerate}
\item  $\mathbf{e}_{i}\le_{s}m\mathbf{e}_{i}$ for all $m\mathbf{e}_{i}\in\mathcal{M}_{d}$,
\item  $m_{1}\mathbf{e}_{i}\le_{s}m_{2}\mathbf{e}_{i}$ implies $m_{1}\cdot s\mathbf{e}_{i}\le_s m_{2}\cdot s\mathbf{e}_{i}$
for all $s\in\mathcal{M},m_{1}\mathbf{e}_{i},m_{2}\mathbf{e}_{i}\in\mathcal{M}_{d}$.
\end{enumerate}
\end{defn}

For convenience, $=_s$ is replaced by $=$. In fact, the admissible order $\le_{s}$ implies the following properties.

\begin{prop}\label{well-order} The admissible module order $\le_{s}$
is a well-order on $\mathcal{M}_{d}$, and it extends the order $|$
on $\mathcal{M}_{d}$, i.e., $m_{1}\mathbf{e}_{i}\,|\, m_{2}\mathbf{e}_{i}$ implies
$m_{1}\mathbf{e}_{i}\le_{s}m_{2}\mathbf{e}_{i}$, for all $m_{1}\mathbf{e}_{i},m_{2}\mathbf{e}_{i}\in\mathcal{M}_{d},i\in\{1,\ldots,d\}$.
\end{prop}

\begin{pf}
If $m_{1}\mathbf{e}_{i}\,|\, m_{2}\mathbf{e}_{i}$
in $\mathcal{M}_{d}$, then there exists $t\in\mathcal{M}$ with $t\cdot m_{1}\mathbf{e}_{i}=_{s}m_{2}\mathbf{e}_{i}$.
Since $\mathbf{e}_{i}\le_{s}m_{3}\mathbf{e}_{i}$, this implies
\begin{equation*}
m_{1}\mathbf{e}_{i}=1\cdot m_{1}\mathbf{e}_{i}\le_s t\cdot m_{1}\mathbf{e}_{i}=m_{2}\mathbf{e}_{i}.
\end{equation*}
 This shows that $\le_{s}$ extends $|$ on $\mathcal{M}_{d}$. By
Dickson's lemma, $\le_{s}$ is a Dickson partial order on $\mathcal{M}_{d}$.
And $\le_{s}$ is a well-order on $\mathcal{M}$ as it is a linear
order.
\end{pf}

It should be noticed that $\le_{s}$ may or may not be related to
$\le_{m}$. The \textbf{compatible} property \citet{Kreuzer} between $\le_{m}$
and $\le_{s}$ is used for the proof of termination for the GVWHS
algorithm in \citet{Volny}: $\sigma\mathbf{e}_{j}\le_{s}\tau\mathbf{e}_{j}$ if and only if $\sigma\le_{m}\tau$. And in \citet{Arri11}, this property is implicitly used in the proof of termination. The following section will show that this relation is indispensable for the proof of finiteness.

For any $sp=(\mathbf{u},p)\in\mathcal{SP}$, let $lm_{\le_{s}}(\mathbf{u})=\mu\mathbf{e}_{k}$
be the \textbf{signature} of $(\mathbf{u},p)$ and $lm_{\le_{m}}(p)$
the leading monomial of $(\mathbf{u},p)$. By an abuse of notation, we
write $lm$ for $lm_{\leq_{s}}$ and $lm_{\leq_{m}}$ if no misunderstanding
occurs. We call $k=idx(\mathbf{u})=idx(sp)$ the \textbf{index} and
call $\mu$ the monomial of the signature. The set of the signatures of elements in $\mathcal{SP}^{*}=\mathcal{SP}\setminus\{(\mathbf{0},0)\}$ is denoted by $sig(\mathcal{SP^*})$.

\section{Properties of sig-polynomials}
\label{sig-polynomial}

\begin{defn} Define a map
\setlength{\arraycolsep}{0.0em}
\begin{eqnarray*}
\begin{array}{cccc}
LM: & \mathbf{NSP} & \rightarrow & \mathcal{M}_{d}\times\mathcal{M}\\
 & (\mathbf{u},p) & \rightarrow & (\mathbf{s},m)=(lm(\mathbf{u}),lm(p)),
\end{array}
\end{eqnarray*}
\setlength{\arraycolsep}{5pt}
 and three orders $\prec_{m,s}$, $\prec_{s,m}$ and $|_{super}$ on the image
$LM(\mathbf{NSP})$ in the following way:
\setlength{\arraycolsep}{0.0em}
\begin{eqnarray*}
\begin{array}{ccccc}
(\mathbf{s}',m') \prec_{m,s} (\mathbf{s},m) & \Leftrightarrow & \lambda\cdot m'=m &\ and\ & \lambda\cdot \mathbf{s}'<_{s}\mathbf{s},\\
(\mathbf{s}'',m'')\prec_{s,m}(\mathbf{s},m) & \Leftrightarrow & \lambda\cdot \mathbf{s}''=\mathbf{s} &\ and\ & \lambda\cdot m''<_{m}m,\\
(\mathbf{s}^*,m^*) |_{super} (\mathbf{s},m) & \Leftrightarrow & \lambda\cdot m^{*}=m &\ and\ & \lambda\cdot \mathbf{s}^{*}=\mathbf{s},
\end{array}
\end{eqnarray*}
\setlength{\arraycolsep}{5pt}
 where $sp,sp',sp''\in\mathbf{NSP}$, $\lambda\in\mathcal{M}$,
and by the relation $<_{s}$ ($<_{m}$) is meant the strict part of
the associated admissible module order. \end{defn}

Under the map $LM$, the image of a sig-polynomial is called a \textbf{leading pair}. We can generalize
two orders $\prec_{s,m}$ and $|_{super}$ on $LM(\mathcal{SP}^{*})$ by
adding the following definitions.
\setlength{\arraycolsep}{0.0em}
\begin{eqnarray*}
\begin{array}{ccccc}
(\mathbf{s}',0) & \prec_{s,m} & (\mathbf{s},m) & \Leftrightarrow & \mathbf{s}'\,|\, \mathbf{s},\\
(\mathbf{s}^{*},0) & |_{super} & (\mathbf{s}'',0) & \Leftrightarrow & \mathbf{s}^{*}\,|\, \mathbf{s}'',
\end{array}
\end{eqnarray*}
\setlength{\arraycolsep}{5pt}
 where the sig-polynomials above are all in $\mathcal{SP}^{*}$ and
$m\neq0$.

Without confusion, denote $|$ on $LM(\mathcal{SP}^{*})$ instead of $|_{super}$ too. Now,
a special kind of reduction is introduced as follows.

\begin{defn}[Top-Reduction] Let $(\mathbf{u},p)\in\mathcal{SP}^{*}$ be a sig-polynomial and $\mathbf{B}\subseteq\mathcal{SP}^{*}$ a set of
sig-polynomials. $sp$ is called to be \textbf{top-reducible} by $\mathbf{B}$,
if there exists a sig-polynomial $(\mathbf{u}',p')\in\mathbf{B}$ satisfying one
of the three conditions below,
\begin{enumerate}
\item \label{tm-reduction} $LM(\mathbf{u}',p')\prec_{m,s}LM(\mathbf{u},p)$, for $lm(p)\neq0$,
\item \label{ts-reduction} $LM(\mathbf{u}',p')\prec_{s,m}LM(\mathbf{u},p)$, for $lm(p)\neq0$,
\item \label{super} $LM(\mathbf{u}',p')\,|\, LM(\mathbf{u},p)$;
\end{enumerate}
otherwise, $(\mathbf{u},p)$ is \textbf{top-irreducible} by $\mathbf{B}$. Such
a top-reduction is called \textbf{regular}, if item \ref{tm-reduction} or \ref{ts-reduction} is satisfied,
and \textbf{super} otherwise.
\end{defn}

For convenience, we call a regular top-reduction satisfying item \ref{tm-reduction} by a \textbf{tm-reduction} (top monomial reduction) for short, and call a regular top-reduction satisfying item \ref{ts-reduction} by a \textbf{ts-rewriting}\footnote{The term ``ts-rewriting'' has the similar meaning as the ``M-pair'' in \citet{Volny}.} (top signature rewriting). Let $sp_1\in \mathcal{SP}$ be a non-syzygy sig-polynomial. We say that $sp_1$ is tm-reducible by $\mathcal{SP}^*$ if there exists $sp_2\in\mathcal{SP}^*$ such that $sp_1\underset{\mathcal{SP}^*}{\longrightarrow}sp_2$, i.e., $sp_2$ is the tm-reduction result of $sp_1$ by some sig-polynomial in $\mathcal{SP}^*$. $\xrightarrow[\mathcal{SP}^*]{*}$ is the reflexive-transitive closure of $\underset{\mathcal{SP}^*}{\longrightarrow}$.
Let $\mathcal{SG}$ be a subset of $\mathcal{SP}$. $\mathcal{SG}$ is called an \textbf{S-Gr\"{o}bner basis} for
the module $\mathcal{SP}$, if every nonzero sig-polynomial $sp\in\mathcal{SP}$
is top-reducible by $\mathcal{SG}$. This definition is the same as the one in
\citet{GVW}. Hence, by \citet[Prop. 2.2]{GVW}, define the
Gr\"{o}bner basis for the syzygy module of $\mathbf{f}$ by $\mathcal{G}_{0}=\{\mathbf{u}\,|\,(\mathbf{u},p)\in \mathcal{SG},p=0\}$,
and define the Gr\"{o}bner basis for $\mathcal{I}$ by $\mathcal{G}_1=\{p\,|\,(\mathbf{u},p)\in \mathcal{SG},p\neq0\}$.

Certainly, there exist different S-Gr\"{o}bner bases for a polynomial ideal. Before investigating
S-Gr\"{o}bner bases, let us consider the properties of the order
$\prec_{m,s}$ ($\prec_{s,m}$).

\begin{prop} \label{minimal}
\begin{enumerate}
\item The order $\prec_{m,s}$ ($\prec_{s,m}$)
is strictly well-founded partial-order on $LM(\mathbf{NSP})$ ($LM(\mathcal{SP}^{*})$).

\item Let $S_{p}$ be the set of $\prec_{m,s}$-minimal
elements in $LM(\mathbf{NSP})$ and $S_{q}$ the set of $\prec_{s,m}$-minimal elements in $LM(\mathcal{SP}^{*})$, then $S_{q}=S_{p}\oplus S_{syz}$,
where $S_{syz}=\{(\mathbf{s},m)\in S_{q}\,|\, m=0\}$.
\end{enumerate}
\end{prop}

\begin{pf}
\begin{enumerate}
\item It is easy to see that $\prec_{m,s}$ and $\prec_{s,m}$ are
irreflexive, strictly antisymmetric and transitive. Assume for a contradiction
that the sequence $\{(\mathbf{s}_{n},m_{n})\}_{n\in\mathbb{N}}$ is a strictly descending
$\prec_{m,s}$-chain in $LM(\mathbf{NSP})$, then $m_{i}=m_{N}$ when
$i>N$ for some $N\in\mathbb{N}$ since $|$ is well-founded
on $\mathcal{M}$. For any $i>j>N$, we must have $\mathbf{s}_{i}<_{s}\mathbf{s}_{j}<_{s}\mathbf{s}_{N}$.
Thus, $\{\mathbf{s}_{n}\}_{n\geq N}$ form a strictly descending $<_{s}$-chain
in $sig(\mathcal{SP^*})$, whereas the admissible module order $\le_{s}$ is a well-order,
a contradiction. Similarly, $\prec_{s,m}$ is strictly well-founded
because $|$ is well-founded on $\mathcal{M}_{d}$ and the admissible
module order $\le_{m}$ is a well-order even if $0$ is added.

\item For $(\mathbf{s}_1,m_1)\in S_{p}$, assume that there exists a leading pair $(\mathbf{s}_2,m_2)$
in $LM(\mathcal{SP}^{*})$ such that $(\mathbf{s}_2,m_2)\prec_{s,m}(\mathbf{s}_1,m_1)$.
Then there exists a nonzero monomial $m\in\mathcal{M}$ such that
$m \mathbf{s}_2=\mathbf{s}_1$ and $m m_2<_{m}m_1$.
Let $\mathbf{s}_3=lm(\mathbf{s}_1-m\mathbf{s}_2)$ and $m_3=lm(m_1-m m_2)$, then $m_3=m_1$
and $\mathbf{s}_3<_{s}\mathbf{s}_1$. Thus, there is a leading pair $(\mathbf{s}_3,m_3)\in LM(SP^*)$ such that $(\mathbf{s}_3,m_3)\prec_{m,s}(\mathbf{s}_1,m_1)$,
a contradiction. For $(\mathbf{s}_1,m_1)\in S_{q}\setminus S_{syz}$, it can
be proved similarly that $(\mathbf{s}_1,m_1)$ is also in $S_{p}$.
\end{enumerate}
\end{pf}

Below follows a natural corollary.
\begin{cor}
Let $sp$ be a non-syzygy sig-polynomial in $\mathcal{SP}^*$. $sp$ is ts-rewritable by $\mathcal{SP}^*$ if and only if it is tm-reducible by $\mathcal{SP}^*$.
\end{cor}

It can be seen that $\mathbf{ISP}=\{sp\in\mathcal{SP}^{*}\,|\, LM(sp)\in S_{q}\}$
is the set of all sig-polynomials which are not ts-rewritable by $\mathcal{SP}^*$. Super
top-reducing elements further in $\mathbf{ISP}$ results the subset of all top-irreducible sig-polynomials
called the \textbf{top-reduced S-Gr\"{o}bner basis} $\mathcal{TSG}$ for $\mathcal{SP}$.
The signature of a top-irreducible sig-polynomial is defined by
the \textbf{top-irreducible signature} for $\mathcal{SP}$. Besides, by two
\textbf{equivalent sig-polynomials} $sp$ and $sp'$ we mean $sp'\ne sp$ such that
$LM(sp')=LM(sp)$. If we store only one for equivalent sig-polynomials in $\mathcal{TSG}$,
for fixed orders $\le_{m}$ and $\le_{s}$,
the top-reduced S-Gr\"{o}bner basis $\mathcal{TSG}$ is uniquely determined by the module
$\mathcal{SP}$ up to equivalence. Those top-reducible sig-polynomials in $\mathcal{SP}\setminus \mathcal{TSG}$ are also called \textbf{redundant sig-polynomials}.

Since $(\mathcal{M},|)$ and $(\mathcal{M}_{d}, |)$ are Dickson partial ordered sets, by Lemma \ref{direct product}, we have $(\mathcal{M}_{d}\times\mathcal{M}, |^*)$ is also a Dickson partial ordered set of which the order $|^*$ is defined as follows:
\begin{equation*}
(\mathbf{s}_1,m_1)\,|^*\,(\mathbf{s}_2,m_2)\Leftrightarrow \mathbf{s}_1\,|\,\mathbf{s}_2\ and\ m_1\,|\,m_2,
\end{equation*}
where $(\mathbf{s}_1,m_1),(\mathbf{s}_2,m_2)$ are in $\mathcal{M}_{d}\times\mathcal{M}$.

\begin{lem}\label{compatible}
Let $(\mathbf{s}_1,m_1)$ and $(\mathbf{s}_2,m_2)$ be two arbitrary leading pairs in $LM(\mathcal{SP}^*)$ such that $(\mathbf{s}_1,m_1)\,|^*\,(\mathbf{s}_2,m_2)$. If the admissible monomial order $\leq_m$ and the admissible module order $\leq_s$ are compatible, then $(\mathbf{s}_1,m_1)$ and $(\mathbf{s}_2,m_2)$ are comparable with respect to one of the three orders $\preceq_{m,s}$, $\preceq_{s,m}$ and $|$.
\end{lem}

\begin{pf}
Let $s$ and $m$ be two monomials in $\mathcal{M}$ such that $s=\mathbf{s}_2/\mathbf{s}_1$ and $m=m2/m1$. There are three cases as follows.
\begin{enumerate}
  \item If $m = s$, then $(\mathbf{s}_1,m_1)\,|\,(\mathbf{s}_2,m_2)$.
  \item If $s<_m m$, then $sm_1<_m m_2$, and $(\mathbf{s}_1,m_1)\prec_{s,m}(\mathbf{s}_2,m_2)$.
  \item If $m<_m s$, as $\leq_m$ and $\leq_s$ are compatible, $m\mathbf{s}_1<_s s\mathbf{s}_1 = \mathbf{s}_2$, and $(\mathbf{s}_1,m_1)\prec_{m,s}(\mathbf{s}_2,m_2)$.
\end{enumerate}
Therefore, $(\mathbf{s}_1,m_1)$ and $(\mathbf{s}_2,m_2)$ are comparable with respect to one of the three orders $\preceq_{m,s}$, $\preceq_{s,m}$ and $|$.
\end{pf}

The finiteness of the top-reduced S-Gr\"{o}bner basis is due
to the following fact.

\begin{thm}\label{div}
The divisibility
relation $|$ is a Dickson partial order on $LM(\mathcal{SP}^{*})$.
Moreover, the top-reduced S-Gr\"{o}bner basis for $\mathcal{SP}$ is finite.
\end{thm}

\begin{pf} It is straightforward to verify that $|$ is reflexive,
transitive and antisymmetric. Since $|^*$ is a Dickson partial order on $LM(\mathcal{SP}^*)$, the $|^*$-minimal elements in $LM(\mathcal{SP}^*)$ are finite. Because the leading pair of a top-irreducible sig-polynomial is $|^*$-minimal in $LM(\mathcal{SP}^*)$ by Lemma \ref{compatible}. So there are a finite number of top-irreducible sig-polynomials in $\mathcal{TSG}$ up to equivalence.
\end{pf}

It can be seen that the ``compatible'' property is indispensable for the finiteness of the top-reduced S-Gr\"{o}bner basis $\mathcal{TSG}$. Hence in the remaining sections of this paper, we will assume that the admissible monomial order $\leq_m$ and the admissible module order $\leq_s$ are compatible.
Suppose two sig-polynomials $sp_{1}=(\mathbf{u}_{1},p_{1}), $ $sp_{2}=(\mathbf{u}_{2},p_{2})\in \mathbf{NSP}$.
Let
\begin{equation*}
m=lcm(lm(p_{1}),lm(p_{2})),\ m_{1}=\frac{m}{lm(p_{1})},m_{2}=\frac{m}{lm(p_{2})}.
\end{equation*}
If $m_{1}lm(\mathbf{u}_{1})>_{s}m_{2}lm(\mathbf{u}_{2})$, then
\begin{itemize}
\item $cp=m_{1}(\mathbf{u}_{1},p_{1})=(m_{1}\mathbf{u}_{1},m_{1}p_{1})$ is
called a \textbf{J-pair} of $sp_{1}$ and $sp_{2}$;
\item $sp_1$ ($sp_2$) is called the \textbf{first (second) component} of $cp$;
\item $m_{1}$ and $m_{2}$ are called the \textbf{multipliers} of $sp_{\ensuremath{1}}$
and $sp_{2}$.
\end{itemize}

\section{The GVWHS Algorithm}
\label{the GVWHS}

As in \citet{Volny}, it can be deduced $\mathbf{e}_{1},\ldots,\mathbf{e}_{d}$
are top-irreducible signatures. Let $sp=(\mathbf{e}_{i},g)$ be a
sig-polynomial in $\mathbf{NSP}$, where $1\leq i\leq d$. If $LM(sp)$ is not
$\prec_{s,m}$-minimal in $LM(\mathcal{SP}^{*})$, there must exist a sig-polynomial
$sp'=(\mathbf{e}_{i},g')\in \mathcal{SP}^{*}$ whose leading pair is $\prec_{s,m}$-minimal.
As $sp'$ cannot be super top-reduced by $\mathcal{SP}^{*}$, $sp'$ is top-irreducible sig-polynomial
and $\mathbf{e}_{i}$ is top-irreducible signature.

For a signature $\mathbf{s}$, we denote by $\mathcal{SP}_{\leq_s(\mathbf{s})}$ the subset
of sig-polynomials in $\mathcal{SP}$ of which the signatures are smaller
than or equal to $\mathbf{s}$ with respect to the order $\leq_s$, and denote by $\mathcal{SG}_{\leq_s(\mathbf{s})}$ the S-Gr\"{o}bner basis
for $\mathcal{SP}_{\leq_s(\mathbf{s})}$. We have the following theorem which is similar but stronger than \citet[Th. 4.11]{Volny}.

\begin{thm} \label{sign}
Let $\mathbf{s}$ be a signature in $sig(\mathcal{SP}^*)$ such that $\mathbf{s}\neq\mathbf{e}_{i}$
for any $1\leq i\leq d$. $\mathbf{s}$ is top-irreducible if and only if
$\mathbf{s}$ is the signature of a J-pair $cp$ of two non-syzygy
top-irreducible sig-polynomials with smaller signatures and $cp$
is not ts-rewritable by $\mathcal{SG}_{<_s(\mathbf{s})}$.
\end{thm}

\begin{pf}
As $\prec_q$ is well-founded, and $\mathbf{e}_1,\ldots,\mathbf{e}_d$ are top-irreducible,
there exists a top-irreducible sig-polynomial $(\mathbf{u}_k,g_k)$ such that $m_k(\mathbf{u}_k,g_k)$
has signature $\mathbf{s}$ and $m_k(\mathbf{u}_k,g_k)$ is not ts-rewritable by $\mathcal{SG}_{<_s(\mathbf{s})}$,
where $m_k>_m1$. Because $\mathbf{s}$ is a top-irreducible
signature, $LM(m_{k}\mathbf{u}_{k},m_{k}g_{k})$ is not $\prec_{s,m}$-minimal
in $LM(\mathcal{SP}^{*})$. Hence there is a non-syzygy top-irreducible sig-polynomial
$(\mathbf{u}_{j},g_{j})\in\mathcal{SG}_{<_s(\mathbf{s})}$ tm-reducing $m_{k}(\mathbf{u}_{k},g_{k})$.
Denote by $m_{k}'(\mathbf{u}_{k},g_{k})$ the J-pair of $(\mathbf{u}_{k},g_{k})$
and $(\mathbf{u}_{j},g_{j})$, where $m_{k}'\,|\, m_{k}$.

Assume for a contradiction that $m_{k}'$ properly divides $m_{k}$.
Because $m_{k}'(\mathbf{u}_{k},g_{k})$ is not ts-rewritable by $\mathcal{SG}_{<_s(\mathbf{s})}$,
after a sequence of tm-reduction on $m_{k}'(\mathbf{u}_{k},g_{k})$,
we get a tm-irreducible sig-polynomial $(\mathbf{u}_{i},g_{i})$ and
$g_{i}\neq0$. $(\mathbf{u}_{i},g_{i})$ is equivalent to $m_{t}(\mathbf{u}_{t},g_{t})$,
a monomial multiple of some top-irreducible sig-polynomial $(\mathbf{u}_{t},g_{t})\in\mathcal{SG}_{<_s(\mathbf{s})}$,
where $m_{t}\geq_{m}1$. Thus, $(\mathbf{u}_{t},g_{t})$ ts-rewrites
$m_{k}'(\mathbf{u}_{k},g_{k})$ and thence $m_{k}(\mathbf{u}_{k},g_{k})$,
a contradiction. Therefore, $m_{k}'=m_{k}$, that is, $cp=(m_{k}\mathbf{u}_{k},m_{k}g_{k})$
is the J-pair of two non-syzygy top-irreducible sig-polynomials with
smaller signatures such that $\mathbf{s}=lm(m_{k}\mathbf{u}_{k})=lm(\mathbf{u}')$
and $cp$ is not ts-rewritable by $\mathcal{SG}_{<_s(\mathbf{s})}$.

For the forward direction, assume for a contradiction that $\mathbf{s}$
is not top-irreducible. Then $LM(\mathbf{u},g)$ is $\prec_{s,m}$-minimal
in $LM(\mathcal{SP}^{*})$. By Proposition \ref{minimal}, $LM(\mathbf{u},g)$
is also $\prec_{m,s}$-minimal in $LM(\mathcal{SP}^{*})$ as $g\neq0$, a contradiction.
\end{pf}

First, we present the GVWHS algorithm, which is modified slightly
from the algorithm mentioned in \citet{Volny}. The subset of non-syzygy sig-polynomials in $\mathcal{SG}$ is denoted by $\mathbf{G}_1$ and $sig(\mathbf{G}_1)$ is the set of signatures of sig-polynomials in $\mathbf{G}_1$. Let $\mathbf{S}$ be a set of polynomials (sig-polynomials), sort($\mathbf{S}$, $\leq_m$ ($\leq_s$)) means that
we arrange $\mathbf{S}$ by ascending leading monomials (signatures) of polynomials (sig-polynomials) with respect to the order $\leq_m$ ($\leq_s$).

\begin{algorithm}
\caption{The GVWHS algorithm}
\begin{algorithmic}[1]
\State \textbf{inputs:}
    \Statex $\mathbf{F}=\{f_{1},\ldots,f_{d}\}\in R$, a list of polynomials
    \Statex $\le_{m}$ an admissible monomial order on $\mathcal{M}$
    \Statex $\le_{s}$, an admissible module order on $\mathcal{M}_{d}$ which is compatible with $\leq_m$
\State \textbf{outputs:}
\Statex $\mathcal{G}_1$, a Gr\"{o}bner basis for $\mathcal{I}=<f_{1},\ldots,f_{d}>$
\State interreduce $\mathbf{F}$ and $\mathbf{F}:=$sort($\{f_{1},\ldots,f_{d}\}$, $\leq_{m}$)
\State \textbf{init1:}
\Statex $\mathbf{CPs}:=\{(\mathbf{e}_{1},f_{1}),\ldots,(\mathbf{e}_{d},f_{d})\}$
and $\mathcal{SG}=\{(f_{i}\mathbf{e}_{j}-f_{j}\mathbf{e}_{i},0)\,|\,1\le i<j\le d\}$
\While{$\mathbf{CPs}\ne\emptyset$}
     \State \label{min-GVWHS}$cp := $min($\{cp\in \mathbf{CPs}\}$, $\le_{s}$) and $\mathbf{CPs}:=\mathbf{CPs}\backslash\{cp\}$
     \If{$cp$ is not ts-rewritable by $\mathcal{SG}$}
         \State \label{reduce-GVWHS}$cp\xrightarrow[\mathcal{SG}]{*}sp=(\mathbf{u},g)$
         \If{$g\neq 0$}
            \State \label{store-GVWHS}$\mathbf{CPs} := $sort($\mathbf{CPs}\cup\{J-pair(sp,sp')$ $\,|\,\forall sp'\in\mathbf{G}_1\}$, $\le_{s}$) and store only one J-pair for each distinct signature of minimal leading monomial
            \State $\mathcal{SG}:=\mathcal{SG}\cup\{(g\mathbf{u}_{l}-g_l\mathbf{u},0)\,|\, (\mathbf{u}_{l},g_l)\in \mathbf{G}_1\}$
         \EndIf
         \State $\mathcal{SG}:=\mathcal{SG}\cup{sp}$
    \EndIf
\EndWhile
\State \textbf{return} $\{g\,|\, (\mathbf{u},g)\in \mathcal{SG}\setminus \mathbf{Syz}\}$
\end{algorithmic}
\end{algorithm}

The only difference compared with the basic
algorithm in \citet{Volny} is that we discard the Gr\"{o}bner basis for
the syzygy module when the algorithm terminates. Note that
the proof for correctness in \citet{Volny} is not complete. Suppose
$\mathbf{s}$ is a top-irreducible signature, it must be
proved, as in Theorem \ref{sign}, that there exists a J-pair of $cp$ such that $cp$ is an M-pair and $\mathbf{s}_{cp}=\mathbf{s}$.

\begin{thm} For any finite subset $\mathbf{F}$ of polynomials in $R$, the
GVWHS algorithm terminates after finitely many steps and it
creates a Gr\"{o}bner basis for the ideal $\mathcal{I}=<\mathbf{F}>$.
\end{thm}

\begin{pf}
We proceed by induction on the top-irreducible signature $\mathbf{s}$.
Because $<_{s}$ is an admissible module order on $\mathcal{M}_{d}$,
the smallest signature of sig-polynomials in $\mathcal{SP}^{*}$ must be one
of the top-irreducible signatures $\mathbf{e}_{1},\ldots,\mathbf{e}_{d}$,
denoted by $\mathbf{e}_{i}$. The case $\mathbf{s}=\mathbf{e}_{i}$ is trivial.
As $Cps$ is initialized with $\{(\mathbf{e}_{1},f_{1}),\ldots,(\mathbf{e}_{d},f_{d})\}$,
during the first while-loop, $(\mathbf{e}_{i},f_{i})$ is added into
$\mathcal{SG}$, which is the S-Gr\"{o}bner basis $\mathcal{SG}_{\leq_s(\mathbf{e}_i)}$.

Let $\mathbf{s}>\mathbf{e}_{i}$, and suppose that $\mathcal{SG}$ created by the GVWHS
algorithm is $\mathcal{SG}_{<_s(\mathbf{s})}$ after finitely many while-loops. If $\mathbf{s}=\mathbf{e}_{j}$, where $1\leq j\leq d$, $j\neq i$, there
exists a $cp=(\mathbf{e}_{j},f_{j})$ at line \ref{min-GVWHS} and $cp$ is not ts-rewritable
by $\mathcal{SG}_{<_s(\mathbf{e}_j)}$. Tm-reducing $cp$ repeatedly by $\mathcal{SG}_{<_s(\mathbf{e}_j)}$
at line \ref{reduce-GVWHS} results a top-irreducible sig-polynomial $sp$ with signature
$\mathbf{e}_{j}$ because $\mathbf{e}_{j}$ is top-irreducible. Thus,
$\mathcal{SG}_{\leq_s(\mathbf{e}_j)}$ can be obtained. If $\mathbf{s}\neq\mathbf{e}_{j}$,
we can also obtain a J-pair $cp'$ with signature $\mathbf{s}$ at line \ref{min-GVWHS}
and $cp'$ is not ts-rewritable by $\mathcal{SG}_{<_s(\mathbf{s})}$ by Theorem \ref{sign}.
After that, a top-irreducible sig-polynomial $sp$ with signature
$\mathbf{s}$ is created and $\mathcal{SG} = \mathcal{SG}_{\leq_s(\mathbf{s})}$. Because top-irreducible signatures are finite in
$\mathcal{SP}$, after finitely many steps, $\mathcal{SG}=\mathcal{SG}_{\leq_s(\mathbf{s}_{max})}$ is the S-Gr\"{o}bner basis.

By Theorem \ref{sign}, the remaining J-pairs
in $\mathbf{CPs}$, if any, are all sig-polynomials with top-irreducible signatures and they will be ts-rewritten by $\mathcal{SG}$.
Therefore, the algorithm terminates and generates $\mathcal{SG}$,
an S-Gr\"{o}bner basis for $\mathcal{SP}$, and the output is a Gr\"{o}bner basis for the ideal
$\mathcal{I}=<\mathbf{F}>$.
\end{pf}

In the remaining part of this section, we aim to reformulate the GVWHS algorithm
into an F5G algorithm (F5-like algorithm in GVWHS's style) and find out the connection between the GVWHS
algorithm and the F5 algorithm. It is, we shall see, an
F5-like algorithm with a different insertion stategy. Before proceeding
to prove the termination of the F5G algorithm, we introduce another order as follows.
%we define an order like \citet{Huang} between two sig-polynomials.

\begin{defn} Define an order $\preceq_l$ on $LM(\mathcal{SP}^*)$ in the following way:
\begin{equation*}
    (\mu\mathbf{e}_i,m)\preceq_l (\mu'\mathbf{e}_i,m') \Leftrightarrow m\mu'\leq_m m'\mu
\end{equation*}
\end{defn}
Note that the order $\preceq_l$ is not defined when two elements in $LM(\mathcal{SP}^*)$ are with different signatures, so $\preceq_{l}$ is a well-founded quasi-order on $LM(\mathcal{SP}^*)$. Particularly, if we restrict $\preceq_l$ on the subset $\{(\mu\mathbf{e}_i,m)\in LM(\mathcal{SP}^*)\,|\,i=i_0,1\leq i_0 \leq d)\}$, then $\preceq_l$ is a well-order on it. Moreover, if $(\mu\mathbf{e}_i,m)\prec_{s,m} (\mu'\mathbf{e}_i,m')$ or $(\mu'\mathbf{e}_i,m')\prec_{m,s} (\mu\mathbf{e}_i,m)$, we have $(\mu\mathbf{e}_i,m)\prec_l (\mu'\mathbf{e}_i,m')$.

Below is the pseudo code of the F5G algorithm in Buchberger's style
which is similar to the algorithm in \citet{SunAMSS09}. We detach the set $\mathbf{Psyz}$ of principal syzygies from $\mathcal{SG}$, and the remainder is denoted by $\mathcal{SG}'$. That is to say, $\mathcal{SG}=\mathbf{PSyz}\cup\mathcal{SG}'$. As is known that there may exist syzygies in $\mathcal{SG}'$, so by $\mathbf{G}_1$ is meant the set of non-syzygy sig-polynomials in $\mathcal{SG}'$. The notations are similar with those in the GVWHS algorithm.

\begin{algorithm}
\caption{The F5G Algorithm (F5-like algorithm in GVWHS's style)}
\begin{algorithmic}[1]
\State \textbf{inputs:}
    \Statex $\mathbf{F}=\{f_{1},\ldots,f_{d}\}\in R$, a list of polynomials
    \Statex $\le_{m}$ an admissible monomial order on $\mathcal{M}$
    \Statex $\le_{s}$, an admissible module order on $\mathcal{M}_{d}$ which is compatible with $\leq_m$
    \Statex $\preceq_{l}$, an order on $LM(\mathcal{SP}^*)$
\State \textbf{outputs:}
\Statex $\mathcal{G}_1$, a Gr\"{o}bner basis for $\mathcal{I}=<f_{1},\ldots,f_{d}>$
\State interreduce $\mathbf{F}$ and $\mathbf{F}:=$sort($\{f_{1},\ldots,f_{d}\}$, $\leq_{m}$),
$F_{i}=(\mathbf{e}_{i,}f_{i})$ for $i=1,\ldots,d$
\State \textbf{init2:}
\Statex $\mathbf{CPs}:=$sort($\{J-pair[F_{i},F_{j}]\,|\,1\leq i<j\leq d\}$, $\leq_s$), $\mathcal{SG}'=\{F_{i}\,|\, i=1,\ldots,d\}$ and
$\mathbf{PSyz}=\{(f_{i}\mathbf{e}_{j}-f_{j}\mathbf{e}_{i},0)\,|\,1\le i<j\le d\}$
\While{$\mathbf{CPs}\ne\emptyset$}
     \State \label{min-F5G}$cp:=$min($\{cp\in \mathbf{CPs}\}$, $\le_{s}$) and $\mathbf{CPs}:=\mathbf{CPs}\backslash\{cp\}$
     \If{$cp$ is neither ts-rewritable by $\mathbf{PSyz}$ nor F5-rewritable
by $\mathcal{SG}'$}\label{criteria-F5G}
         \State \label{reduce-F5G}$cp\xrightarrow[\mathcal{SG}']{*}sp=(\mathbf{u},g)$
         \State $\mathcal{SG}':=$insert\_by\_decreasing\_l($sp$, $\mathcal{SG}'$, $\preceq_{l}$)
         \If{$g\neq 0$}
            \State \label{store-F5G}$\mathbf{CPs}:=$sort($\mathbf{CPs}\cup\{J-pair(sp,sp')$ $\,|\,\forall sp'\in \mathbf{G}_1, sp'\neq sp\}$, $\leq_{s}$)$=\{mSG'(k)\}$ and store only one J-pair for each distinct signature of which the first component has maximum index $k$ in $\mathcal{SG}'$
            \State $\mathbf{PSyz}:=\mathbf{PSyz}\cup\{(g\mathbf{u}_{l}-g_l\mathbf{u},0)$ $\,|\, (\mathbf{u}_{l},g_l)\in \mathbf{G}_1\}$ and discard those super top-reducible in $\mathbf{PSyz}$
         \EndIf
    \EndIf
\EndWhile
\State \textbf{return} $\{g\,|\, (\mathbf{u},g)\in \mathcal{SG}'\setminus \mathbf{Syz}\}$
\end{algorithmic}
\end{algorithm}

\algnewcommand\algorithmicto{\textbf{to}}
\algrenewtext{For}[2]%
    {\algorithmicfor\ #1 \algorithmicto\ #2 \algorithmicdo}
\begin{algorithm}
\caption{F5-rewritable}
\begin{algorithmic}[1]
\State \textbf{inputs:}
\Statex $cp=m(\mathbf{u}_{k},g_{k})\in \mathcal{SP}$
\Statex $\mathcal{SG}':=\mathcal{SG}'(i)=\{(\mathbf{u}_{1},g_{1}),\ldots,(\mathbf{u}_{r},g_{r})\}$
\State \textbf{outputs:}
\Statex \textbf{true} if $m\mathbf{u}_k$ is F5-rewritable by another sig-polynomial in $\mathcal{SG}'$
\State find the first index $j_{b}$ and the last index $j_{e}$ in
$\mathcal{SG}'$ such that $idx(sp)=idx(\mathcal{SG}'(j_{b}))=idx(\mathcal{SG}'(j_{e}))$
\For{$i=j_{e}$}{$j_{b}$}
    \If{$lm(\mathbf{u}_{i})\,|\, lm(m\mathbf{u}_{k})$}
        \State \textbf{return} $i\neq k$
    \EndIf
\EndFor
\State \textbf{return} \textbf{false}
\end{algorithmic}
\end{algorithm}

\begin{algorithm}
\caption{insert\_by\_decreasing\_l}
\begin{algorithmic}[1]
\State \textbf{inputs:}
\Statex $sp$, a sig-polynomial
\Statex $\mathcal{SG}':=\mathcal{SG}'(i)=\{(\mathbf{u}_{1},g_{1}),\ldots,(\mathbf{u}_{r},g_{r})\}$
\Statex $\preceq_l$, an order on $LM(\mathcal{SP}^*)$
\State find the first index $j_{b}$ and the last index $j_{e}$ in
$\mathcal{SG}'$ such that $idx(sp)=idx(\mathcal{SG}'(j_{b}))=idx(\mathcal{SG}'(j_{e}))$
\For{$i=j_{e}$}{$j_{b}$}
    \If{$LM(\mathcal{SG}'(i))\succeq_{l}LM(sp)$}
        \State insert $sp$ into $\mathcal{SG}'$ after $\mathcal{SG}'(i)$
        \State \textbf{return}
    \EndIf
\EndFor
\State insert $sp$ into $\mathcal{SG}'$ before $\mathcal{SG}'(j_b)$
\State \textbf{return}
\end{algorithmic}
\end{algorithm}

It is important to note that the index $k$ mentioned at line \ref{store-F5G}, different from the index of a sig-polynomial, points to the sig-polynomial of the $k$-th position in $\mathcal{SG}'$.

Let $sp_j$, $sp_i$ be two sig-polynomials in $\mathcal{SG}'$ and let $cp=tsp_i$, $cp'=t'sp_j$ be two J-pairs with the same signature. From the insert\_by\_decreasing\_l function, we know that $sp_j$ appears later in $\mathcal{SG}'$ than $sp_i$ if $LM(sp_j)\prec_l LM(sp_i)$. In this case, $cp'$ is discarded as its first component $sp_j$ is ahead of the first component $sp_i$ of $cp$. Hence line \ref{store-F5G} of the F5G algorithm is equivalent to storing only one J-pair for each distinct signature of minimal leading monomial at line \ref{store-GVWHS} of the GVWHS algorithm. Even more, the F5G algorithm adopts the same criterion as the GVWHS algorithm for finding redundant sig-polynomials.

\begin{lem}\label{ts-rewritable}
During an execution of the while-loop,
let $cp_{0}$ be the J-pair chosen at line \ref{min-F5G} in the F5G algorithm, and let $\mathbf{CPs}_0$ be the value of $\mathbf{CPs}$, $\mathbf{PSyz}_{0}$ the value
of $\mathbf{PSyz}$, and $\mathcal{SG}'_{0}$ the value of $\mathcal{SG}'$ at line \ref{min-F5G}. The criteria
of line \ref{criteria-F5G} in the F5G algorithm are equivalent to the statement of
judging whether $cp_{0}$ is not ts-rewritable by $\mathbf{PSyz}_{0}\cup \mathcal{SG}'_{0}$.
\end{lem}

\begin{pf}
Assume that the J-pair $cp_{0}=m(\mathbf{u}_{k},g_{k})=msp_{k}$
is ts-rewritable by $\mathcal{SG}'_{0}$ in the F5G algorithm. We may find $sp_{i}=(\mathbf{u}_{i},g_{i})\in \mathcal{SG}'_{0}$
ts-rewrite $cp_{0}$ and $lm(m_{i}\mathbf{u}_{i})=lm(m\mathbf{u}_{k})$,
where $m_{i}>_{m}1$. Since $LM(sp_{i})\prec_{l}LM(sp_{k})$ means $i>k$,
$cp_{0}$ is F5-rewritable by $sp_i$ in the F5G algorithm. That is to say, $cp_0$ can not pass the criteria of line \ref{criteria-F5G} in the F5G algorithm.

If $cp_{0}=msp_{k}\in \mathbf{CPs}_0$ is not ts-rewritable by $\mathcal{SG}'_{0}$. Assume for
a contradiction that $cp_{0}$ is F5-rewritable by $sp_{j}\in \mathcal{SG}'_{0}$,
$j>k$. We know $LM(sp_{j})\neq_{l}LM(sp_{k})$, or else the J-pair $cp_{0}$ had been discarded
by line\ref{store-F5G} of the F5G algorithm. So $LM(sp_{j})\prec_{l}LM(sp_{k})$, which means
$lm(m_{j}\mathbf{u}_{j})=lm(m\mathbf{u}_{k})$ and $lm(m_{j}g_{j})<lm(mg_{k})$,
where $m_{j}>_{m}1$. Hence $cp_{0}$ is ts-rewritable by $sp_{j}$,
a contradiction.
\end{pf}

Note that Lemma \ref{ts-rewritable} does not apply to the algorithms we will discuss later since the insertion strategy of the F5G is used for the proof. In Theorem \ref{sign}, two components of a J-pair have to be
top-irreducible. As a matter of fact, a generalized lemma follows.

\begin{lem} \label{gen-sign} If $\mathbf{s}$ is the signature of a J-pair
$cp=msp_{k}=m(\mathbf{u}_{k},g_{k})$ of two non-syzygy sig-polynomials
$sp_{k}$ and $sp_{j}$ (with smaller signatures) and $cp$ is not
ts-rewritable by $\mathcal{SG}_{<_s(\mathbf{s})}$, then $\mathbf{s}$ is a top-irreducible
signature of $\mathcal{SP}$.
\end{lem}

\begin{pf}
Assume for a contradiction that $lm(m\mathbf{u}_{k})$
is not a top-irreducible signature. Then $LM(cp)$ is $\prec_{s,m}$-minimal
in $LM(\mathcal{SP}^{*})$. But there exists $m'sp_{j}=m(\mathbf{u}_{j},g_{j})$
such that $lm(m'\mathbf{u}_{j})<_{s}lm(m\mathbf{u}_{k})$ and $lm(m'g_{j})=lm(mg_{k})$,
that is, $cp$ is tm-reducible by $\mathcal{SG}_{<_s(\mathbf{s})}$, a contradiction.
\end{pf}

\begin{thm}
For any finite subset $\mathbf{F}$ of polynomials in $R$, the
F5G algorithm terminates after finitely many steps and it creates
a Gr\"{o}bner basis for the ideal $\mathcal{I}=<\mathbf{F}>$.
\end{thm}

\begin{pf} Due to Lemma \ref{ts-rewritable}, we will use the criterion of judging whether $cp$ is not ts-rewritable by $\mathbf{PSyz}\cup\mathcal{SG}'$ instead.
Similar to the corresponding proof of the GVWHS algorithm,
we proceed by induction on the top-irreducible signature $\mathbf{s}$. Because
$<_{s}$ is an admissible module order on $\mathcal{M}_{d}$, the
smallest signature of sig-polynomials in $\mathcal{SP}^{*}$ must be one of
the top-irreducible signatures $\mathbf{e}_{1},\ldots,\mathbf{e}_{d}$,
denoted by $\mathbf{e}_{i}$. The case $\mathbf{s}=\mathbf{e}_{i}$ is trivial.
As $\mathcal{SG}'$ is initialized with $\{(\mathbf{e}_{1},f_{1}),\ldots,(\mathbf{e}_{d},f_{d})\}$,
$\mathcal{SG}'$ is the S-Gr\"{o}bner basis for $\mathcal{SP}_{\leq_s(\mathbf{e}_i)}$.

Let $\mathbf{s}>\mathbf{e}_{i}$, and suppose that $\mathbf{PSyz}\cup \mathcal{SG}'$ created by
the F5G algorithm is $\mathbf{PSyz}_{<_s(\mathbf{e}_j)}\cup\mathcal{SG}'_{<_s(\mathbf{e}_j)}=\mathcal{SG}_{<_s(\mathbf{s})}$
after finitely many while-loops. If $\mathbf{s}=\mathbf{e}_{j}$, where $1\leq j\leq d$, $j\neq i$, there
is only one sig-polynomial $(\mathbf{e}_{j},f_{j})$ in $\mathcal{SG}'_{<_s(\mathbf{e}_j)}$ with top-irreducible
signature $\mathbf{e}_{j}$. And if
$(\mathbf{e}_{j},f_{j})$ is tm-irreducible by $\mathcal{SG}_{<_s(\mathbf{e}_j)}$,
$\mathbf{PSyz}_{<_s(\mathbf{e}_j)}\cup\mathcal{SG}'_{<_s(\mathbf{e}_j)}$ is $\mathcal{SG}_{\leq_s(\mathbf{e}_j)}$.
If $(\mathbf{e}_{j},f_{j})$ is tm-reducible by $\mathcal{SG}_{<_s(\mathbf{e}_j)}$,
during an execution of the while-loop, line \ref{min-F5G} will create a J-pair
$cp=(\mathbf{e}_{j},f_{j})$ and $cp$ is not ts-rewritable by $\mathcal{SG}_{<_s(\mathbf{e}_j)}$.
Tm-reducing $cp$ repeatedly by $\mathcal{SG}_{<_s(\mathbf{e}_j)}$ at line
\ref{reduce-F5G} results a top-irreducible sig-polynomial $sp$ with signature $\mathbf{e}_{j}$
because $\mathbf{e}_{j}$ is top-irreducible. Thus, $\mathcal{SG}_{\leq_s(\mathbf{e}_j)}$
can be obtained. If $\mathbf{s}\neq\mathbf{e}_{j}$, we can also obtain a
J-pair $cp'$ with signature $\mathbf{s}$ at line \ref{min-F5G} and $cp'$ is not ts-rewritable
by $\mathcal{SG}_{<_s(\mathbf{s})}$ by Theorem \ref{sign}. After that, a top-irreducible
sig-polynomial $sp$ with signature $\mathbf{s}$ will be created. Because
top-irreducible signatures are finite in $\mathcal{SP}$, after finitely many
steps, $\mathbf{PSyz}\cup\mathcal{SG}'=\mathcal{SG}_{\leq_s(\mathbf{s}_{max})}$ is the S-Gr\"{o}bner basis.

By Lemma \ref{gen-sign}, the remaining J-pairs in $\mathbf{CPs}$, if any, are all ts-rewritable by $\mathcal{SG}$.
Therefore, the algorithm terminates and generates an S-Gr\"{o}bner basis $\mathbf{PSyz}\cup \mathcal{SG}'$
for $\mathcal{SP}$, and the output is a Gr\"{o}bner basis for the ideal $\mathcal{I}=<\mathbf{F}>$.
\end{pf}

\section{The termination and correctness of the F5B Algorithm}
\label{the F5B}

We present two variants of the F5 algorithm here and in the next section, both of which share the same
F5-rewritten criterion with that in the F5G algorithm. So we do not write the F5-rewritable function in detail again.

For two non-syzygy components $sp_{1}$ and $sp_{2}$ of a J-pair,
let $m_{1}$ and $m_{2}$, respectively, be their multipliers. A much
simpler version than the F5B algorithm (F5 algorithm in Buchberger's style) in \citet{SunAMSS09} is given below. The F5B
algorithm here does not apply F5-rewritable check for $m_{2}sp_{2}$
nor in the tm-reduction of the J-pair. Omitting these influences neither
the termination nor the correctness of the F5B algorithm in \citet{SunAMSS09}.
For details, one can refer to \citet{SunISSAC11} and \citet{EderISSAC11}.

\begin{algorithm}
\caption{The F5B Algorithm (F5 algorithm in Buchberger's style)}
\begin{algorithmic}[1]
\State \textbf{inputs:}
    \Statex $\mathbf{F}=\{f_{1},\ldots,f_{d}\}\in R$, a list of polynomials
    \Statex $\le_{m}$ an admissible monomial order on $\mathcal{M}$
    \Statex $\le_{s}$, an admissible module order on $\mathcal{M}_{d}$ which is compatible with $\leq_m$
    \Statex $\preceq_{l}$, an order on $LM(\mathcal{SP}^*)$
\State \textbf{outputs:}
\Statex $\mathcal{G}_1$, a Gr\"{o}bner basis for $\mathcal{I}=<f_{1},\ldots,f_{d}>$
\State interreduce $\mathbf{F}$ and $\mathbf{F}:=$sort($\{f_{1},\ldots,f_{d}\}$, $\leq_{m}$),
$F_{i}=(\mathbf{e}_{i,}f_{i})$ for $i=1,\ldots,d$
\State \textbf{init2:}
\Statex $\mathbf{CPs}:=$sort($\{J-pair[F_{i},F_{j}]\,|\,1\leq i<j\leq d\}$, $\leq_s$), $\mathcal{SG}'=\{F_{i}\,|\, i=1,\ldots,d\}$ and
$\mathbf{PSyz}=\{(f_{i}\mathbf{e}_{j}-f_{j}\mathbf{e}_{i},0)\,|\,1\le i<j\le d\}$
\While{$\mathbf{CPs}\ne\emptyset$}
     \State \label{min-F5B}$cp:=$min($\{cp\in \mathbf{CPs}\}$, $\le_{s}$) and $\mathbf{CPs}:=\mathbf{CPs}\backslash\{cp\}$
     \If{$cp$ is neither ts-rewritable by $\mathbf{PSyz}$ nor F5-rewritable
by $\mathcal{SG}'$}\label{criteria-F5B}
         \State \label{reduce-F5B}$cp\xrightarrow[\mathcal{SG}']{*}sp=(\mathbf{u},g)$
         \State $\mathcal{SG}':=$insert\_by\_index($sp$, $\mathcal{SG}'$)
         \If{$g\neq 0$}
            \State \label{store-F5B}$\mathbf{CPs}:=$sort($\mathbf{CPs}\cup\{J-pair(sp,sp')$ $\,|\,\forall sp'\in \mathbf{G}_1, sp'\neq sp\}$, $\leq_{s}$)$=\{mSG'(k)\}$ and store only one J-pair for each distinct signature of which the first component has maximum index $k$ in $\mathcal{SG}'$
            \State $\mathbf{PSyz}:=\mathbf{PSyz}\cup\{(g\mathbf{u}_{l}-g_l\mathbf{u},0)$ $\,|\, (\mathbf{u}_{l},g_l)\in \mathbf{G}_1\}$ and discard those super top-reducible in $\mathbf{PSyz}$
         \EndIf
    \EndIf
\EndWhile
\State \textbf{return} $\{g\,|\, (\mathbf{u},g)\in \mathcal{SG}'\setminus \mathbf{Syz}\}$
\end{algorithmic}
\end{algorithm}

\begin{algorithm}
\caption{insert\_by\_index}
\begin{algorithmic}[1]
\State \textbf{inputs:}
\Statex $sp$, a sig-polynomial
\Statex $\mathcal{SG}':=\mathcal{SG}'(i)=\{(\mathbf{u}_{1},g_{1}),\ldots,(\mathbf{u}_{r},g_{r})\}$

\State find the last index $j_{e}$ in
$\mathcal{SG}'$ such that $idx(sp)=idx(\mathcal{SG}'(j_{e}))$
\State insert $sp$ into $\mathcal{SG}'$ after $\mathcal{SG}'(j_{e})$
\State \textbf{return}
\end{algorithmic}
\end{algorithm}

Instead of using an auxiliary number for each sig-polynomial in \citet{SunAMSS09}, the F5B algorithm here realizes the same rewritable check by adjusting the order of sig-polynomials in $\mathcal{SG}'$. One can find that the real difference between the F5B and F5G algorithms is the insertion of elements in $\mathcal{SG}'$.
The reason why line \ref{store-F5B} does not affect the correctness of the algorithm lies in the fact that the first component of the discarded J-pair appears earlier in $\mathcal{SG}'$ than that of the stored J-pair.

\begin{lem} \label{more-gen-sign}Let $\mathbf{s}$ be a signature in $sig(\mathcal{SP}^*)$ such that $\mathbf{s}\neq\mathbf{e}_{i}$ for any $1\leq i\leq d$.
During an execution of the while-loop in the F5B algorithm, let $\mathbf{PSyz}_{<_s(\mathbf{s})}$ and $\mathcal{SG}'_{<_s(\mathbf{s})}$
be the values of $\mathbf{PSyz}$ and $\mathcal{SG}'$. If $\mathbf{s}$ is top-irreducible, then $\mathbf{s}$ is the signature
of a J-pair $cp$ of two non-syzygy sig-polynomials in $\mathcal{SG}'_{<_s(\mathbf{s})}$ with smaller signatures
and $cp$ is neither ts-rewritable by $\mathbf{PSyz}_{<_s(\mathbf{s})}$ nor F5-rewritable
by $\mathcal{SG}'_{<_s(\mathbf{s})}$.
\end{lem}

\begin{pf}
By Theorem \ref{sign}, there exists a J-pair $cp'=m'(\mathbf{u}_{k},g_{k})$
of two non-syzygy top-irreducible sig-polynomials with smaller signatures
such that $\mathbf{s}=lm(m'\mathbf{u}_{k})$ and $cp'$ is not ts-rewritable
by $\mathcal{SG}_{<_s(\mathbf{s})}$. If $cp'$ is not F5-rewritable by $\mathcal{SG}'_{<_s(\mathbf{s})}$,
$cp'$ is the desired $cp$. If $cp'$ is F5-rewritable by $\mathcal{SG}'_{<_s(\mathbf{s})}$,
let $(\mathbf{u}_{j},g_{j})$ be the non-syzygy sig-polynomial in
$\mathcal{SG}'_{<_s(\mathbf{s})}$ F5-rewriting $cp'$ as $\mathbf{s}$ is top-irreducible
signature. That is, $(\mathbf{u}_{j},g_{j})$ satisfies that $lm(m\mathbf{u}_{j})=\mathbf{s}$
and $m(\mathbf{u}_{j},g_{j})$ is ts-rewritable by $\mathcal{SG}_{<_s(\mathbf{s})}$.
Further more, $(\mathbf{u}_{j},g_{j})$ is the sig-polynomial in $\mathcal{SG}'_{<_s(\mathbf{s})}$ with
the largest signature dividing $\mathbf{s}$ according to the structure of
the F5B algorithm. As $\mathbf{s}$ is top-irreducible signature, $m(\mathbf{u}_{j},g_{j})$
is not ts-rewritable by the principal syzygy submodule $\mathbf{PSyz}_{<_s(\mathbf{s})}$
and it can be tm-reduced by some non-syzygy top-irreducible sig-polynomial
$(\mathbf{u}_{t},g_{t})$ in $\mathcal{SG}'_{<_s(\mathbf{s})}$. Denote by $m^{*}(\mathbf{u}_{j},g_{j})$
the J-pair of $(\mathbf{u}_{j},g_{j})$ and $(\mathbf{u}_{t},g_{t})$,
where $m^{*}\,|\, m$.

Assume for a contradiction that $m^{*}$ properly divides $m$. It
can be deduced that $m^{*}(\mathbf{u}_{j},g_{j})$ is neither ts-rewritable
by $\mathbf{PSyz}_{<_s(\mathbf{s})}$ nor F5-rewritable by $\mathcal{SG}'_{<_s(\mathbf{s})}$. After a
sequence of tm-reduction on $m^{*}(\mathbf{u}_{j},g_{j})$, we get
a tm-irreducible sig-polynomial $(\mathbf{u}_{v},g_{v})$ added later
in $\mathcal{SG}'_{<_s(\mathbf{s})}$ than $(\mathbf{u}_{j},g_{j})$. Because $\mathbf{u}_{v}\,|\, m^{*}\mathbf{u}_{j}\,|\, \mathbf{s}$
and $lm(\mathbf{u}_{j})\leq_{s}lm(\mathbf{u}_{v})$, which contradict
the fact that $(\mathbf{u}_{j},g_{j})$ F5-rewrites $cp'$. Therefore,
$m^{*}=m$, that is, $cp=(m\mathbf{u}_{j},mg_{j})$ is the J-pair
of two non-syzygy sig-polynomials with smaller signatures such that
$\mathbf{s}=lm(m\mathbf{u}_{j})=lm(\mathbf{u}')$ and $cp$ is neither ts-rewritable
by $\mathbf{PSyz}_{<_s(\mathbf{s})}$ nor F5-rewritable by $\mathcal{SG}'_{<_s(\mathbf{s})}$.
\end{pf}

It is important to note, however, that we can not guarantee the reverse direction of Lemma \ref{more-gen-sign} is satisfied too. That is to say, there may exist a J-pair $cp$ such that it passes the criteria and the signature of $cp$ is top-reducible. This situation does exist by running experiments: tm-reducing $cp$ will result a redundant sig-polynomial $sp$ which is super top-reducible another computed sig-polynomial.
But the order $\preceq_l$ can be employed for the proof of the termination of the F5B algorithm. Assume the algorithm has created the S-Gr\"{o}bner basis $\mathcal{SG}$
after finite while-loops. Let $sp_i$ and $sp_j$ be two sig-polynomials in $\mathcal{SG}'$ such that $i<j$ i.e., $sp_i$ appears earlier
in $\mathcal{SG}'$ than $sp_j$. We call $(sp_i,sp_j)$ a \textbf{misplaced pair} if $LM(sp_i)\prec_l LM(sp_j)$.
Note that we always order $sp_i$ before $sp_j$ in the misplaced pair. Clearly, the misplacement is the reason for the J-pair of the form  $msp_j$.
\begin{defn} Let $(sp_i,sp_j)$ and $(sp_k,sp_l)$ be two misplaced pairs. And define $(sp_i,sp_j)\prec_{pm}(sp_k,sp_l)$, if one of the following cases is satisfied.
\begin{enumerate}
\item $LM(sp_i)\prec_lLM(sp_k)$
\item $LM(sp_i)\prec_lLM(sp_k)$ and $LM(sp_j)\prec_lLM(sp_l)$
\end{enumerate}
\end{defn}
If each J-pair $msp_j$ is either ts-rewritable by $\mathbf{PSyz}$ or F5-rewritable
by $\mathcal{SG}'$, we call the misplaced pair $(sp_i,sp_j)$ is \textbf{corrected}.

\begin{thm}\label{proof-F5B} For any finite subset $\mathbf{F}$ of polynomials in $R$, the
F5B algorithm terminates after finitely many steps and it creates
a Gr\"{o}bner basis for the ideal $\mathcal{I}=<\mathbf{F}>$.
\end{thm}

\begin{pf}
We still proceed by induction on the top-irreducible signature
$\mathbf{s}$. If $\mathbf{s}=\mathbf{e}_i$ is the smallest top-irreducible signature, the initialized $\mathbf{PSyz}\cup\mathcal{SG}'$ is
the S-Gr\"{o}bner basis $\mathcal{SG}_{\leq_s(\mathbf{e}_i)}$.

Let $\mathbf{s}>\mathbf{e}_{i}$, and suppose that $\mathbf{PSyz}\cup \mathcal{SG}'$ created by
the F5B algorithm is $\mathbf{PSyz}_{<_s(\mathbf{s})}\cup\mathcal{SG}'_{<_s(\mathbf{s})}=\mathcal{SG}_{<_s(\mathbf{s})}$ after
finitely many while-loops. If $\mathbf{s}=\mathbf{e}_{j}$, where $1\leq j\leq d$, $j\neq i$, there
is only one sig-polynomial $(\mathbf{e}_{j},f_{j})$ in $\mathcal{SG}'_{<_s(\mathbf{e}_j)}$ with top-irreducible
signature $\mathbf{e}_{j}$. And if
$(\mathbf{e}_{j},f_{j})$ is tm-irreducible by $\mathcal{SG}_{<_s(\mathbf{e}_j)}$,
$\mathbf{PSyz}_{<_s(\mathbf{e}_j)}\cup\mathcal{SG}'_{<_s(\mathbf{e}_j)}$ is $\mathcal{SG}_{\leq_s(\mathbf{e}_j)}$.
If $(\mathbf{e}_{j},f_{j})$ is tm-reducible by $\mathcal{SG}_{<_s(\mathbf{e}_j)}$,
during an execution of the while-loop, line \ref{min-F5B} will create a J-pair
$cp=(\mathbf{e}_{j},f_{j})$ and $cp$ is neither ts-rewritable by
$\mathbf{PSyz}_{<_s(\mathbf{e}_j)}$ nor F5-rewritable by $\mathcal{SG}'_{<_s(\mathbf{e}_j)}$.
Tm-reducing $cp$ repeatedly by $\mathcal{SG}_{<_s(\mathbf{e}_j)}$ at line
\ref{reduce-F5B} results a top-irreducible sig-polynomial $sp$ with signature $\mathbf{e}_{j}$
because $\mathbf{e}_{j}$ is top-irreducible. Thus, $\mathcal{SG}_{\leq_s(\mathbf{e}_j)}$
can be obtained. If $\mathbf{s}\neq\mathbf{e}_{j}$, we can also obtain a
J-pair $cp'$ with signature $\mathbf{s}$ at line \ref{min-F5B} and $cp'$ is neither
ts-rewritable by $\mathbf{PSyz}_{<_s(\mathbf{s})}$ nor F5-rewritable by $\mathcal{SG}'_{<_s(\mathbf{s})}$
by Lemma \ref{more-gen-sign}. After that, a top-irreducible sig-polynomial
$sp$ with signature $\mathbf{s}$ will be created. Because top-irreducible
signatures are finite in $\mathcal{SP}$, after finitely many steps, the algorithm
generates an S-Gr\"{o}bner basis $\mathbf{PSyz}\cup \mathcal{SG}'=\mathcal{SG}$ for $\mathcal{SP}$.

If there are J-pairs in $\mathbf{CPs}$ at this time, a new $cp''=m(\mathbf{u}_k,g_k)=msp_k$ may pass the criteria and thus generating
a new tm-irreducible sig-polynomial $sp_n$ in $\mathcal{SG}'$. There must exist a top-irreducible $sp_h$ in $\mathcal{SG}'$ such that $sp_h$ can super top-reduce $sp_n$ and
$(sp_h,sp_k)$ is a misplaced pair. That is, $LM(sp_h)=_lLM(sp_n)\prec_lLM(sp_k)$ and $h<k<n$. On the one hand, the J-pairs of $sp_n$ and other possible sig-polynomials, be of the form $m'sp_n$ or not, will generate tm-irreducible
sig-polynomials, say, $sp_p$ with $\preceq_l$-smaller leading pairs if it passes the criteria
of the F5B algorithm. Since the leading pair of $sp_p$ is equal to that of a top-irreducible sig-polynomial and the
top-irreducible sig-polynomials in $\mathcal{SP}$ are finite, this process of creating a J-pair and generating a sig-polynomial always terminates. On the other hand, after finite steps, the misplaced pair $(sp_h,sp_k)$
will be corrected. Though an insertion of a new tm-irreducible sig-polynomial may produce other misplaced pairs, the $\prec_{pm}$-maximum misplaced pair of $\mathcal{SG}'$ without being corrected
is gradually decreasing with respect to the order $\prec_{pm}$. As there are finite pairs not $\prec_{pm}$-equal, the algorithm will terminate finally and output a
Gr\"{o}bner basis for $\mathcal{I} = <\mathbf{F}>$.
\end{pf}

\section{Proof of the termination of the F5 algorithm}
\label{the F5}

In the original F5 algorithm in \citet{Fau02}, the input polynomials
in $F=\{f_{1},\ldots,f_{d}\}$ are homogeneous, and after initialization,
sig-polynomials are $(\mathbf{e}_{1},f_{1}),\ldots,(\mathbf{e}_{d},f_{d})$.
A property follows: If $sp=(\mathbf{u},g)\in \mathbf{NSP}$ and $idx(\mathbf{u})=i$,
$1\leq i\leq d$, then $deg(lm(\mathbf{u}))+deg(lm(f_{i}))=deg(g)$.
We define the \textbf{g-weighted degree} the same with that in \citet{GVW}:
The g-weighted degree $gw-deg$ of a sig-polynomial $sp=(\mathbf{u},g)$
is equal to $deg(lm(\mathbf{u}))+deg(lm(f_{idx(\mathbf{u})}))$. Therefore,
selecting critical pairs of the minimal degree in the original F5
algorithm equals selecting J-pairs of the minimal g-weighted degree.
For an admissible monomial order $\leq_m$, we define the admissible module order $\leq_{s_{0}}$ as
follows.

We say that $x^{\alpha}\mathbf{e}_{i}<_{s_0}x^{\beta}\mathbf{e}_{j}$
if
\begin{enumerate}
\item  $i<j$,
\item  $i=j$ and $gw-deg(x^{\alpha}\mathbf{e}_{i})<gw-deg(x^{\beta}\mathbf{e}_{j})$,
\item  $i=j$ , $gw-deg(x^{\alpha}\mathbf{e}_{i})=gw-deg(x^{\beta}\mathbf{e}_{j})$
and $x^{\alpha}<_{m}x^{\beta}$.
\end{enumerate}
Particularly, we have $x^{\alpha}\mathbf{e}_{i}=_{s_0}x^{\beta}\mathbf{e}_{j}$,
if $i=j$ and $x^{\alpha}=x^{\beta}$.

Sure enough, the order $\leq_{s_{0}}$ is an admissible module order.
By using this order $\leq_{s_{0}}$, we can understand the reformulation
of the original F5 algorithm easier.
In \citet{Fau02}, Faug\`{e}re builded up an array Rule to store the ordering of sig-polynomials on which the F5-rewritten criterion is based. As presented in the following pseudo code, we will just discard the Rule and store the ordering directly in $\mathcal{SG}'$.

Though the F5B and original F5 algorithms share the same F5-rewritten criterion, the ordering in $\mathbf{SG}'$ of the F5B algorithm slightly differs from that in Rule of the F5 algorithm. In the F5B algorithm, let $sp_{1}$ and $sp_{2}$ be two sig-polynomials of the same index
in $\mathcal{SG}'$. If $\mathbf{s}_{sp_{1}}<_{s}\mathbf{s}_{sp_{2}}$,
$sp_{1}$ must appear earlier in $\mathcal{SG}'$ than $sp_{2}$. This is also
interpreted as an isRewritten criterion in \citet{Hashemi11}. However,
in the original F5 algorithm, the claim is not true for sig-polynomials.
Since the Rule is updated not only in the Spol function of \citet{Fau02} but also in the TopReduction function, at the end of each run though the while-loop, the newly added sig-polynomials in Rule have the same index. Moreover, the g-weighted degrees of them are equal as the input polynomials of the original F5 algorithm are homogeneous. Then there is no guarantee that the sig-polynomials are arranged in $\leq_{s_0}$-descending order (note that the original F5 algorithm insert new sig-polynomials at the beginning of Rule). By running several examples, this non-monotony in Rule is verified.

Nevertheless, a weaker relation exists between sig-polynomials
in Rule. During an execution of the while-loop in the original F5
algorithm, let $d$ be the minimal degree of critical pairs. The sig-polynomials
added in Rule are all of g-weighted degree $d$ in the Spol and TopReduction functions.
Hence if two sig-polynomials $sp_{1}$
and $sp_{2}$ in Rule are of the same index satisfying $gw-deg(sp_{1})<gw-deg(sp_{2})$,
then $sp_{1}$ appears earlier in Rule than $sp_{2}$. Besides, if
a J-pair $cp$ of two non-syzygy sig-polynomials $sp_{3}$ and $sp_{4}$
passes criteria of the original F5 algorithm and it is F5-reduced\footnote{Here, F5-reducing means using F5-rewritable check and tm-reducing.} to
$sp_{5}$, then $sp_{5}$ appears later than $sp_{3}$ and $sp_{4}$.
Here the latter property plays an important part in the proof below.

The following is the F5GEN algorithm (F5 algorithm with a generalized insertion strategy) derived from the original one
in \citet{Fau02}. It use the same F5-rewritten criterion as the previous ones. Here we still omit F5-rewritable check when tm-reducing J-pairs as in \citet{Fau02}.

\begin{algorithm}
\caption{The F5GEN Algorithm (F5 algorithm with a generalized insertion strategy)}
\begin{algorithmic}[1]
\State \textbf{inputs:}
    \Statex $\mathbf{F}=\{f_{1},\ldots,f_{d}\}\in R$, a list of polynomials
    \Statex $\le_{m}$ an admissible monomial order on $\mathcal{M}$
    \Statex $\le_{s}$, an admissible module order on $\mathcal{M}_{d}$ which is compatible with $\leq_m$
    \Statex $\preceq_{l}$, an order on $LM(\mathcal{SP}^*)$
\State \textbf{outputs:}
\Statex $\mathcal{G}_1$, a Gr\"{o}bner basis for $\mathcal{I}=<f_{1},\ldots,f_{d}>$
\State interreduce $\mathbf{F}$ and $\mathbf{F}:=$sort($\{f_{1},\ldots,f_{d}\}$, $\leq_{m}$),
$F_{i}=(\mathbf{e}_{i,}f_{i})$ for $i=1,\ldots,d$
\State \textbf{init2:}
\Statex $\mathbf{CPs}:=$sort($\{J-pair[F_{i},F_{j}]\,|\,1\leq i<j\leq d\}$, $\leq_s$), $\mathcal{SG}'=\{F_{i}\,|\, i=1,\ldots,d\}$ and
$\mathbf{PSyz}=\{(f_{i}\mathbf{e}_{j}-f_{j}\mathbf{e}_{i},0)\,|\,1\le i<j\le d\}$
\While{$\mathbf{CPs}\ne\emptyset$}
     \State \label{min-F5GEN}$cp:=$min($\{cp\in \mathbf{CPs}\}$, $\le_{s}$) and $\mathbf{CPs}:=\mathbf{CPs}\backslash\{cp\}$
     \If{$cp$ is neither ts-rewritable by $\mathbf{PSyz}$ nor F5-rewritable
by $\mathcal{SG}'$}\label{criteria-F5GEN}
         \State \label{reduce-F5GEN}$cp\xrightarrow[\mathcal{SG}']{*}sp=(\mathbf{u},g)$
         \State $\mathcal{SG}':=$insert\_F5GEN($sp$, $\mathcal{SG}'$, $cp$)
         \If{$g\neq 0$}
            \State \label{store-F5GEN}$\mathbf{CPs}:=$sort($\mathbf{CPs}\cup\{J-pair(sp,sp')$ $\,|\,\forall sp'\in \mathbf{G}_1, sp'\neq sp\}$, $\leq_{s}$)$=\{mSG'(k)\}$ and store only one J-pair for each distinct signature of which the first component has maximum index $k$ in $\mathcal{SG}'$
            \State $\mathbf{PSyz}:=\mathbf{PSyz}\cup\{(g\mathbf{u}_{l}-g_l\mathbf{u},0)$ $\,|\, (\mathbf{u}_{l},g_l)\in \mathbf{G}_1\}$ and discard those super top-reducible in $\mathbf{PSyz}$
         \EndIf
    \EndIf
\EndWhile
\State \textbf{return} $\{g\,|\, (\mathbf{u},g)\in \mathcal{SG}'\setminus \mathbf{Syz}\}$
\end{algorithmic}
\end{algorithm}

\begin{algorithm}
\caption{insert\_F5GEN}
\begin{algorithmic}[1]
\State \textbf{inputs:}
\Statex $sp$, a sig-polynomial
\Statex $\mathcal{SG}':=\mathcal{SG}'(i)=\{(\mathbf{u}_{1},g_{1}),\ldots,(\mathbf{u}_{r},g_{r})\}$
\Statex $cp=m(\mathbf{u}_k,g_k)$, the J-pair which is tm-reduced to $sp$
\State find the first index $j_{b}$ and the last index $j_{e}$ in
$\mathcal{SG}'$ such that $idx(sp)=idx(\mathcal{SG}'(j_{b}))=idx(\mathcal{SG}'(j_{e}))$
\State insert $sp$ into $\mathcal{SG}'$ after $\mathcal{SG}'(i)$ , where $j_{b}-1\leq i\leq j_{e}$, such that $sp$ appears later in $\mathcal{SG}'$ than $sp_k=(\mathbf{u}_k,g_k)$
\State \textbf{return}
\end{algorithmic}
\end{algorithm}

In the insert\_F5GEN function of the F5GEN algorithm, we can restrict an
appropriate strategy of insertion such that the ordering in $\mathcal{SG}'$
is the same as that in Rule of the original F5 algorithm. The idea for constructing signature-based algorithms also for non-homogeneous polynomial ideals
has been mentioned in Eder and Perry's earlier papers. We shall see that
this F5GEN algorithm here is true for any polynomial ideals both homogeneous and non-homogeneous,
admissible module orders other than $\leq_{s_{0}}$ and the weak condition of ordering in $\mathcal{SG}'$ mentioned in the above pseudo code. But once
the input polynomials are homogeneous and the admissible module order $\leq_{s_{0}}$
is chose, the F5GEN algorithm with an appropriate strategy of insertion
will simulate the original F5 algorithm accurately. Together with the
analysis of equivalence between the original F5 algorithm and the
F5B algorithm in \citet{SunAMSS09}, the proof of termination and correctness
for this F5GEN algorithm can be used to prove the termination and correctness
of the original algorithm in \citet{Fau02}.

\begin{lem} \label{more-gen-sign-1} Let $\mathbf{s}$ be a signature in $sig(\mathcal{SP}^*)$ such that $\mathbf{s}\neq\mathbf{e}_{i}$, for any $1\leq i\leq d$. During an
execution of the while-loop in the F5GEN algorithm, let $\mathbf{PSyz}_{<_s(\mathbf{s})}$
and $\mathcal{SG}'_{<_s(\mathbf{s})}$ be the values of $\mathbf{PSyz}$ and $\mathcal{SG}'$. If $\mathbf{s}$ is top-irreducible, then $\mathbf{s}$
is the signature of a J-pair $cp$ of two non-syzygy sig-polynomials in $\mathcal{SG}'_{<_s(\mathbf{s})}$ with smaller signatures and $cp$ is neither ts-rewritable
by $\mathbf{PSyz}_{<_s(\mathbf{s})}$ nor F5-rewritable by $\mathcal{SG}'_{<_s(\mathbf{s})}$.

\end{lem}

\begin{pf} By Theorem \ref{sign}, there exists a J-pair $cp'=m'(\mathbf{u}_{k},g_{k})$
of two non-syzygy top-irreducible sig-polynomials with smaller signatures
such that $\mathbf{s}=lm(m'\mathbf{u}_{k})$ and $cp'$ is not ts-rewritable
by $\mathcal{SG}_{<_s(\mathbf{s})}$. If $cp'$ is F5-rewritable by $\mathcal{SG}'_{<_s(\mathbf{s})}$,
let $(\mathbf{u}_{j},g_{j})$ be the  last non-syzygy sig-polynomial in
$\mathcal{SG}'_{<_s(\mathbf{s})}$ with the signature dividing $\mathbf{s}$ according to the structure of the
F5GEN algorithm. As $\mathbf{s}$ is top-irreducible signature, $m(\mathbf{u}_{j},g_{j})$
can be tm-reduced by some non-syzygy top-irreducible sig-polynomial
$(\mathbf{u}_{t},g_{t})$ in $\mathcal{SG}'_{<_s(\mathbf{s})}$. Denote by $m^{*}(\mathbf{u}_{j},g_{j})$
the J-pair of $(\mathbf{u}_{j},g_{j})$ and $(\mathbf{u}_{t},g_{t})$,
where $m^{*}\,|\, m$.

Assume for a contradiction that $m^{*}$ properly divides $m$. It
can be deduced that $m^{*}(\mathbf{u}_{j},g_{j})$ is neither ts-rewritable
by $\mathbf{PSyz}_{<_s(\mathbf{s})}$ nor F5-rewritable by $\mathcal{SG}'_{<_s(\mathbf{s})}$. After a
sequence of tm-reduction on $m^{*}(\mathbf{u}_{j},g_{j})$, we get
a tm-irreducible sig-polynomial $(\mathbf{u}_{v},g_{v})$. Because
of the insertion strategy of the F5GEN algorithm, $(\mathbf{u}_{v},g_{v})$ must appear later in $\mathcal{SG}'_{<_s(\mathbf{s})}$ than
$(\mathbf{u}_{j},g_{j})$, which contradicts the fact that $(\mathbf{u}_{j},g_{j})$
F5-rewrites $cp'$. Therefore, $m^{*}=m$, that is, $cp=(m\mathbf{u}_{j},mg_{j})$
is the J-pair of two non-syzygy sig-polynomials with
smaller signatures such that $\mathbf{s}=lm(m\mathbf{u}_{j})=lm(\mathbf{u}')$
and $cp$ is neither ts-rewritable by $\mathbf{PSyz}_{<_s(\mathbf{s})}$ nor F5-rewritable
by $\mathcal{SG}'_{<_s(\mathbf{s})}$.
\end{pf}

For the original F5 algorithm, Gash \citet{Gash} made a conjecture that there is not a sig-polynomial in $\mathcal{SG}'$ super top-reducible by another one. But this can not be satisfied sometimes.
Here we can not guarantee the reverse direction of Lemma \ref{more-gen-sign-1} is satisfied too. It is highly possible that there exist a misplaced pair $(sp_i,sp_j)$ in $\mathcal{SG}'$ as the insertion strategy of the F5 algorithm (it can be seen as a implementation of the F5GEN algorithm) is different from the F5G. From the proof of Theorem \ref{proof-F5B}, we know that $LM(sp)\prec_l LM(sp')\prec_l LM(sp'')$ if the sig-polynomial $sp$ is the result tm-reduced from a J-pair of $sp'$ and $sp''$. A J-pair of the form $msp_j$ may pass the criteria and be reduced to $m'sp_i$ since $LM(sp_i)\prec_l LM(sp_j)$. The sig-polynomial $sp_i$ is added earlier in $\mathcal{SG}'$ than $sp_j$ and $sp_i$ can not be selected in the F5-rewritable function. So both $sp_i$ and $m'sp_i$ are kept in $\mathcal{SG}'$, a contradiction. One can verify this situation by running several examples.

\begin{thm} For any finite subset $\mathbf{F}$ of polynomials in $R$, the
F5GEN algorithm terminates after finitely many steps and it creates
a Gr\"{o}bner basis for the ideal $\mathcal{I}=<\mathbf{F}>$.
\end{thm}

\begin{pf}
Again, we proceed by induction on the top-irreducible
signature $\mathbf{s}$ and let $\mathbf{e}_i$ be the smallest top-irreducible signature.
The case $\mathbf{s}=\mathbf{e}_{i}$ is trivial.

Let $\mathbf{s}>\mathbf{e}_{i}$, and suppose that $\mathbf{PSyz}\cup \mathcal{SG}'$ created by
the F5GEN algorithm is $\mathcal{SG}_{<_s(\mathbf{s})}$ after
finitely many while-loops. If $\mathbf{s}=\mathbf{e}_{j}$, $\mathcal{SG}_{\leq_s(\mathbf{e}_j)}$
can be obtained in similar fashion with the proof of the F5B algorithm. If $\mathbf{s}\neq\mathbf{e}_{j}$, we can also obtain a J-pair $cp'$ with signature $\mathbf{s}$ at line \ref{min-F5GEN} and $cp'$ is neither
ts-rewritable by $\mathbf{PSyz}_{<_s(\mathbf{s})}$ nor F5-rewritable by $\mathcal{SG}'_{<_s(\mathbf{s})}$
by Lemma \ref{more-gen-sign-1}. After that, a top-irreducible sig-polynomial
$sp$ with signature $\mathbf{s}$ will be created. Because top-irreducible
signatures are finite in $\mathcal{SP}$, after finitely many steps, the algorithm
generates an S-Gr\"{o}bner basis $\mathbf{PSyz}\cup \mathcal{SG}'=\mathcal{SG}$ for $\mathcal{SP}$.

If there are J-pairs in $\mathbf{CPs}$ at this time, the leading pair of a newly generated sig-polynomial which is tm-reduced from the J-pair $cp$, is $\preceq_l$-smaller than two components of $cp$. On the one hand, by the insertion strategy of the F5GEN algorithm and leading pair of generated sig-polynomials are $\preceq_l$-equal to that of top-irreducible sig-polynomials, a branch of creating a J-pair and generating a sig-polynomial will end finitely.
On the other hand, after finite steps, a misplaced pair will be corrected. Though an insertion of a new tm-irreducible sig-polynomial may produce other misplaced pairs, the $\prec_{pm}$-maximum misplaced pair of $\mathcal{SG}'$ without being corrected
is gradually decreasing with respect to the order $\prec_{pm}$. As there are finite pairs not $\prec_{pm}$-equal, the algorithm will terminate finally and output a
Gr\"{o}bner basis for $\mathcal{I} = <\mathbf{F}>$.
\end{pf}

Therefore, for any finite set of homogeneous polynomials, the original
F5 algorithm in \citet{Fau02} terminates finitely and it creates
a Gr\"{o}bner basis for the polynomial ideal.

\section{Conclusion}
\label{conclusion}

This paper present a clear proof of the termination of the GVWHS, F5B and F5 algorithms under the condition that the admissible monomial order and the admissible module order are compatible. Of course, there exist some optimizations for improving the efficiency, like recording $(lm(\mathbf{u}), g)$ for each $(\mathbf{u}, g)$ in the implementation. These optimizations do not affect the correctness and termination.
One may find out that the F5G ,F5B and original F5 algorithms are implementations of the F5GEN algorithm with different insertion strategy. That means, the GVWHS algorithm is just an F5-like algorithm. Moreover, with this proved F5GEN algorithm, researchers can shift their focus on the different variants of the F5GEN algorithm and find out the fastest one.

\section*{Acknowledgment}

We would like to thank Christian Eder for valuable feedback and discussions on this work. We would also like to thank Yao Sun, Dingkang Wang and Dongxiao Ma whose comments greatly improved this paper.
\bibliographystyle{elsart-harv}
\bibliography{IEEE_termination}

\end{document}